\documentclass[]{amsart}

\usepackage[initials]{amsrefs}
\usepackage{amsthm}
\usepackage{combelow}
\usepackage{mathtools}
\usepackage{paralist}

\newtheorem{lem}{Lemma}[section]
\newtheorem{prop}[lem]{Proposition}
\newtheorem{theo}[lem]{Theorem}
\newtheorem{cor}[lem]{Corollary}

\theoremstyle{definition}
\newtheorem{definition}[lem]{Definition}

\theoremstyle{remark}
\newtheorem*{rem}{Remark}

\numberwithin{equation}{section}

\providecommand{\norm}[1]{\left\lVert#1\right\rVert}
\providecommand\restr[2]{\left.#1\right|_{#2}}
\renewcommand{\Re}{\operatorname{Re}}
\DeclareMathOperator{\lin}{lin}
\DeclareMathOperator{\spec}{spec}
\DeclareMathOperator{\supp}{supp}
\DeclareMathOperator{\co}{co}
\DeclareMathOperator{\van}{van}

\begin{document}

\title{The Daugavet property and translation-invariant subspaces}
\author[S. L\"ucking]{Simon L\"ucking}
\address{Department of Mathematics\\ Freie Universit\"at Berlin\\ 
         Arnimallee 6\\ 14195 Berlin\\ Germany}
\email{simon.luecking@fu-berlin.de}

\begin{abstract}
    Let $G$ be an infinite, compact abelian group and let $\varLambda$ be a 
    subset of its dual group $\varGamma$. We study the question which spaces of
    the form $C_\varLambda(G)$ or $L^1_\varLambda(G)$ and which quotients of the
    form $C(G)/C_\varLambda(G)$ or $L^1(G)/L^1_\varLambda(G)$ have the Daugavet 
    property.
    
    We show that $C_\varLambda(G)$ is a rich subspace of $C(G)$ if and only if
    $\varGamma \setminus \varLambda^{-1}$ is a semi-Riesz set. If 
    $L^1_\varLambda(G)$ is a rich subspace of $L^1(G)$, then $C_\varLambda(G)$ 
    is a rich subspace of $C(G)$ as well. Concerning quotients, we prove that 
    $C(G)/C_\varLambda(G)$ has the Daugavet property, if $\varLambda$ is a 
    Rosenthal set, and that $L^1_\varLambda(G)$ is a poor subspace of $L^1(G)$, 
    if $\varLambda$ is a nicely placed Riesz set.
\end{abstract}

\subjclass[2010]{Primary 46B04; secondary 43A46}
\keywords{Daugavet property; rich subspace; poor subspace; semi-Riesz set;
          nicely placed set; translation-invariant subspace}
\maketitle

\section{Introduction}

I.\,K.~Daugavet \cite{Daugavet} proved in 1963 that all compact operators $T$
on $C[0,1]$ fulfill the norm identity
\begin{equation*}
    \norm{\mathrm{Id} + T} = 1 + \norm{T},
\end{equation*}
which has become known as the \emph{Daugavet equation}. C.~Foia{\cb{s} and
I.~Singer \cite{FoiasSingerPointsDiffusion} extended this result to all 
weakly compact operators on $C[0,1]$ and A.~Pe{\l}czy{\'n}ski 
\cite{FoiasSingerPointsDiffusion}*{p.~446} observed that their argument can
also be used for weakly compact operators on $C(K)$ provided that $K$ is a
compact space without isolated points. 
Shortly afterwards, G.\,Ya.~Lozanovski{\u\i}
\cite{LozanovskiiAlmostIntegralOperators} showed that the Daugavet equation
holds for all compact operators on $L^1[0,1]$ and J.\,R.~Holub
\cite{HolubDaugavetsEquationL1} extended this result to all weakly compact
operators on $L^1(\varOmega,\varSigma,\mu)$ where $\mu$ is a $\sigma$-finite 
non-atomic measure. V.\,M.~Kadets, R.\,V.~Shvidkoy, G.\,G.~Sirotkin, and 
D.~Werner \cite{KadetsShvidkoySirotkinWernerDaugavetProperty} proved that the 
validity of the Daugavet equation for weakly compact operators already follows 
from the corresponding statement for operators of rank one. This result led to 
the following definition: A Banach space $X$ is said to have the \emph{Daugavet
property}, if every operator $T: X \rightarrow X$ of rank one satisfies the
Daugavet equation. 

Examples include the aforementioned spaces $C(K)$ and 
$L^1(\varOmega,\varSigma,\mu)$, certain function algebras such as the disk 
algebra $A(\mathbb D)$ or the algebra of bounded analytic functions $H^\infty$
\citelist{\cite{WernerDaugavetEquationFunctionSpaces}
\cite{WojtaszczykRemarksDaugavetEquation}}, and non-atomic $C^*$-algebras
\cite{OikhbergDaugavetPropertyC*Algebren}. If $X$ has the Daugavet property,
not only all weakly compact operators on $X$ satisfy the Daugavet equation 
but also all strong Radon-Nikod\'ym operators
\cite{KadetsShvidkoySirotkinWernerDaugavetProperty}, meaning operators $T$ for
which $\overline{T[B_X]}$ is a Radon-Nikod\'ym set, and operators not fixing a
copy of $\ell^1$ \cite{ShvydkoyGeometricAspectsDaugavetProperty}. Furthermore,
$X$ fails the Radon-Nikod\'ym property
\cite{WojtaszczykRemarksDaugavetEquation}, contains a copy of $\ell^1$
\cite{KadetsShvidkoySirotkinWernerDaugavetProperty}, does not have an
unconditional basis \cite{KadetsRemarksDaugavetEquation}, and does not even
embed into a space with an unconditional basis 
\cite{KadetsShvidkoySirotkinWernerDaugavetProperty}.

The listed properties give the impression that spaces with the Daugavet 
property are ``big''. It is therefore an interesting question which subspaces
of a space $X$ with the Daugavet property inherit this property. One approach
is to look at closed subspaces $Y$ such that the quotient space $X/Y$ is
``small''. For this purpose, V.\,M.~Kadets and M.\,M.~Popov 
\cite{KadetsPopovDaugavetPropertyNarrowOperatorsRichSubspaces} introduced
on $C[0,1]$ and $L^1[0,1]$ the class of \emph{narrow} operators, a
generalization of the class of compact operators, and called a subspace 
\emph{rich}, if the corresponding quotient map is narrow. This concept was
transferred to spaces with the Daugavet property by V.\,M.~Kadets, 
R.\,V.~Shvidkoy, and D.~Werner \cite{KadetsShvidkoyWernerNarrowOperators}. Rich 
subspaces inherit the Daugavet property and the class of narrow operators
includes all weakly compact operators, all strong Radon-Nikod\'ym operators,
and all operators which do not fix copies of $\ell^1$ 
\cite{KadetsShvidkoyWernerNarrowOperators}.

If $Y$ is a rich subspace of a Banach space $X$ with the Daugavet property,
then not only $Y$ inherits the Daugavet property but also every closed subspace
of $X$ which contains $Y$. In view of this property, V.\,M.~Kadets, 
V.~Shepelska, and D.~Werner introduced a similar notion for quotients of $X$ and
called a closed subspace $Y$ \emph{poor}, if $X/Z$ has the Daugavet property for
every closed subspace $Z \subset Y$. They also showed that poverty is a dual 
property to richness \cite{KadetsShepelskaWernerQuotientsDaugavetProperty}.

Let us consider an infinite, compact abelian group $G$ with its Haar measure
$m$. Since $G$ has no isolated points and since $m$ has no atoms, the spaces
$C(G)$ and $L^1(G)$ have the Daugavet property. Using the group structure of
$G$, we can translate functions that are defined on $G$ and look at closed,
translation-invariant subspaces of $C(G)$ or $L^1(G)$. These subspaces can be
described via subsets $\varLambda$ of the dual group $\varGamma$ and are of the
form
\begin{gather*}
    C_\varLambda(G) = \{ f \in C(G): \spec f \subset \varLambda \} 
        \quad \text{and} \quad
        L^1_\varLambda(G) = \{ f \in L^1(G): \spec f \subset \varLambda \},
    \shortintertext{where}
    \spec f = \left\{ \gamma \in \varGamma: \hat f(\gamma) \neq 0 \right\}.
\end{gather*}

We are going to study the question which closed, translation-invariant
subspaces of $C(G)$ and $L^1(G)$ and which quotients of
the form $C(G)/C_\varLambda(G)$ or $L^1(G)/L^1_\varLambda(G)$ have the Daugavet 
property. We will characterize rich, translation-invariant subspaces of 
$C(G)$ and will show that $C_\varLambda(G)$ is rich in $C(G)$, if 
$L^1_\varLambda(G)$ is rich in $L^1(G)$. We will prove that 
$C(G)/C_\varLambda(G)$ has the Daugavet property, if 
$L^1_{\varGamma \setminus \varLambda^{-1}}(G)$ is a rich subspace of $L^1(G)$, 
and that $L^1(G)/L^1_\varLambda(G)$ has the Daugavet property, if 
$C_{\varGamma \setminus \varLambda^{-1}}(G)$ is a rich subspace of $C(G)$. 
We will furthermore identify a big class of poor, translation-invariant 
subspaces of $L^1(G)$.

\section{Preliminaries}

Let $\mathbb T$ be the \emph{circle group}, i.e., the multiplicative group of 
all complex numbers with absolute value one. In the sequel, $G$ will be an
infinite, compact abelian group with addition as group operation
and $e_G$ as identity element. $\mathcal B(G)$ will denote its Borel
$\sigma$-algebra, $m$ its normalized Haar measure, $\varGamma$ its (discrete) 
\emph{dual group}, i.e., the group of all continuous homomorphisms from $G$ into
$\mathbb T$, and $\varLambda$ a subset of $\varGamma$. Linear combinations of
elements of $\varGamma$ are called \emph{trigonometric polynomials} and
we set $T(G) = \lin\varGamma$. We will write $\mathbf{1}_{G}$ for the identity 
element of $\varGamma$, which coincides with the function identically equal to 
one.
\begin{lem}
    \label{LemmaCovering}
    If $O$ is an open neighborhood of $e_G$, then there exists a covering of
    $G$ by disjoint Borel sets $B_1, \dotsc, B_n$ with $B_k - B_k \subset O$
    for $k = 1, \dotsc, n$.
    \begin{proof}
        Let $V$ be an open neighborhood of $e_G$ with $V - V \subset O$. Since
        $G$ is compact, we can choose $x_1, \dotsc, x_n \in G$ with
        $G = \bigcup_{k = 1}^n (x_k + V)$. Set $B_1 = x_1 + V$ and
        $B_k = (x_k + V) \setminus \bigcup_{l = 1}^{k-1} B_k$ for $k = 2,
        \dotsc, n$. Then $B_1, \dotsc, B_n$ is a covering of $G$ by disjoint
        Borel sets and for every $k \in \{1, \dotsc, n\}$
        \begin{equation*}
            B_k - B_k \subset (x_k + V) - (x_k + V) \subset V - V \subset O.
            \qedhere
        \end{equation*}
    \end{proof}
\end{lem}
The spaces $L^1(G)$ and $M(G)$ are commutative Banach algebras with respect to 
convolution and $L^1(G)$ is a closed ideal of $M(G)$ 
\cite{RudinFourierAnalysis}*{Theorem~1.1.7, 1.3.2, and 1.3.5}. 
If $\mu \in M(G)$, its \emph{Fourier-Stieltjes transform} is defined by
\begin{equation*}
    \hat \mu(\gamma) = \int_G \overline\gamma \,d\mu
    \quad (\gamma \in \varGamma),
\end{equation*}
and the map $\mu \mapsto \hat \mu$ is injective, multiplicative and continuous
\cite{RudinFourierAnalysis}*{Theorem~1.3.3 and 1.7.3}.
$L^1(G)$ does not have a unit, unless $G$ is discrete. But we
always have an approximate unit 
\cite{HewittRossAbstractHarmonicAnalysisII}*{Remark~VIII.32.33(c) and 
Theorem~VIII.33.12}.
\begin{prop}
    \label{PropositionApproximateUnit}
    There is a net $(v_j)_{j \in J}$ in $L^1(G)$ with the following properties:
    \begin{enumerate}[\upshape(i)]
        \item $\norm{f - f*v_j}_1 \longrightarrow 0$ for every $f \in L^1(G)$;
        \item $\norm{f - f*v_j}_\infty \longrightarrow 0$ 
            for every $f \in C(G)$;
        \item $v_j \geq 0$, $v_j \in T(G)$ and $\hat v_j \geq 0$ 
            for every $j \in J$;
        \item $\norm{v_j}_1 = 1$ for all $j \in J$;
        \item $\hat v_j(\gamma) \longrightarrow 1$ 
            for every $\gamma \in \varGamma$. 
    \end{enumerate}
\end{prop}
If $f:G \rightarrow \mathbb C$ is a function and $x$ an element of $G$, the 
\emph{translate} $f_x$ of $f$ is defined by
\begin{equation*}
    f_x(y) = f(y - x) \quad (y \in G).
\end{equation*}
A subspace $X$ of $L^1(G)$ or $C(G)$ is called \emph{translation-invariant},
if $X$ contains with a function $f$ all possible translates $f_x$. As 
already mentioned in the introduction, all closed, translation-invariant
subspaces of $C(G)$ or $L^1(G)$ are of the form $C_\varLambda(G)$ or
$L^1_\varLambda(G)$
\cite{HewittRossAbstractHarmonicAnalysisII}*{Theorem~IX.38.7}, where
$\varLambda$ is a subset of $\varGamma$. We define analogously 
$T_\varLambda(G)$, $L^\infty_\varLambda(G)$ and $M_\varLambda(G)$. Note that by
Proposition~\ref{PropositionApproximateUnit} the space $T_\varLambda(G)$ is 
$\norm{\,\cdot\,}_\infty$-dense in $C_\varLambda(G)$
and $\norm{\,\cdot\,}_1$-dense in $L^1_\varLambda(G)$.

We will need the following characterization of the Daugavet property
\cite{KadetsShvidkoySirotkinWernerDaugavetProperty}*{Lemma 2.2}.
\begin{lem}
    \label{LemmaCharacterizationSlices}
    Let $X$ be a Banach space. The following assertions are equivalent:
    \begin{enumerate}[\upshape(i)]
        \item $X$ has the Daugavet property.
        \item For every $x \in S_X$, $x^* \in S_{X^*}$, and $\varepsilon > 0$
            there is some $y \in S_X$ such that $\Re x^*(y) \geq 1 -
            \varepsilon$ and $\norm{x + y} \geq 2 - \varepsilon$.
        \item For every $x \in S_X$, $x^* \in S_{X^*}$, and $\varepsilon > 0$
            there is some $y^* \in S_{X^*}$ such that $\Re y^*(x) \geq 1 -
            \varepsilon$ and $\norm{x^* + y^*} \geq 2 - \varepsilon$.
    \end{enumerate}
\end{lem}

\section{Structure-preserving isometries}

The Daugavet property depends crucially on the norm of a space and is
preserved under isometries but in general not under isomorphisms. Considering 
translation-invariant subspaces of $C(G)$ and $L^1(G)$, it would be useful to 
know isometries that map translation-invariant subspaces onto 
translation-invariant subspaces.
\begin{definition}
    Let $G_1$ and $G_2$ be locally compact abelian groups with dual groups
    $\varGamma_1$ and $\varGamma_2$. Let $H : G_1 \rightarrow G_2$ be
    a continuous homomorphism. The \emph{adjoint homomorphism} 
    $H^* : \varGamma_2 \rightarrow \varGamma_1$ is defined by
    \begin{equation*}
        H^*(\gamma) = \gamma \circ H \quad (\gamma \in \varGamma_2). 
    \end{equation*}
\end{definition}
The adjoint homomorphism $H^*$ is continuous 
\cite{HewittRossAbstractHarmonicAnalysisI}*{Theorem~VI.24.38}, \mbox{$H^{**}=H$}
\cite{HewittRossAbstractHarmonicAnalysisI}*{VI.24.41.(a)}, and 
$H^*[\varGamma_2]$ is dense in $\varGamma_1$ if and only if $H$ is one-to-one
\cite{HewittRossAbstractHarmonicAnalysisI}*{VI.24.41.(b)}.
\begin{lem}
    \label{LemmaMeasurePreserving}
    Let $H : G \rightarrow G$ be a continuous and surjective homomorphism.
    Then $H$ is measure-preserving, i.e., each Borel set $B$ of $G$ satisfies
    $m(H^{-1}[B]) = m(B)$.
    \begin{proof}
        Denote by $\mu$ the push-forward of $m$ under $H$. It is easy to see 
        that $\mu$ is regular and $\mu(G) = 1$. Since the Haar measure is
        uniquely determined, it suffices to show that $\mu$ is 
        translation-invariant.
        
        Fix $B \in \mathcal B(G)$ and $x \in G$. $H$ is surjective and thus
        there is $y \in G$ with $H(y) = x$. It is not difficult to check that 
        $H^{-1}[B + H(y)] = H^{-1}[B] + y$. Using this equality, we get
        \begin{equation*}
            \mu(B + x) = m(H^{-1}[B + H(y)]) = m(H^{-1}[B] + y) = 
            m(H^{-1}[B]) = \mu(B). \qedhere
         \end{equation*}
    \end{proof}
\end{lem}
\begin{prop}
    Let $H : \varGamma \rightarrow \varGamma$ be a one-to-one homomorphism and 
    let $\varLambda$ be a subset of $\varGamma$. Then $C_\varLambda(G) \cong
    C_{H[\varLambda]}(G)$ and $L^1_\varLambda(G) \cong L^1_{H[\varLambda]}(G)$.
    \begin{proof}
        If we define $T : C(G) \rightarrow C(G)$ by
        \begin{equation*}
            T(f) = f \circ H^* \quad (f \in C(G)),
        \end{equation*}
        then $T$ is well-defined and an isometry because $H^*$ is continuous
        and surjective. (Note that $H^*[G]$ is compact and therefore closed.)
        For every trigonometric polynomial $f = \sum_{k = 1}^n a_k \gamma_k$ and
        every $x \in G$ we get
        \begin{equation*}
            T(f)(x) = \sum_{k = 1}^n a_k \gamma_k(H^*(x)) =
            \sum_{k = 1} a_k H(\gamma_k)(x).
        \end{equation*}
        Hence $T$ maps for every $\varLambda \subset \varGamma$ the space
        $T_\varLambda(G)$ onto $T_{H[\varLambda]}(G)$ and by density the space
        $C_\varLambda(G)$ onto $C_{H[\varLambda]}(G)$.
        
        Let us look at the same $T$ but now as an operator from $L^1(G)$ into
        itself. It is again an isometry because $H^*$ is measure-preserving
        by Lemma~\ref{LemmaMeasurePreserving}. It still maps for every
        $\varLambda \subset \varGamma$ the space $T_\varLambda(G)$ onto
        $T_{H[\varLambda]}(G)$ and so by density $L^1_\varLambda(G)$ onto
        $L^1_{H[\varLambda]}(G)$.
    \end{proof} 
\end{prop}
\begin{cor}
    \label{CorollaryHomomorphism}
    Let $H : \varGamma \rightarrow \varGamma$ be a one-to-one homomorphism. If
    $C_\varLambda(G)$ has the Daugavet property, then $C_{H[\varLambda]}(G)$ has 
    the Daugavet property as well. Analogously, if $L^1_\varLambda(G)$ has the
    Daugavet property, then $L^1_{H[\varLambda]}(G)$ has the Daugavet property 
    as well.
\end{cor}
Let us give an example. Every one-to-one homomorphism on $\mathbb Z$ is of the 
form $k \mapsto nk$ where $n \neq 0$ is a fixed integer. So 
$C_\varLambda(\mathbb T) \cong C_{n\varLambda}(\mathbb T)$ and 
$L^1_\varLambda(\mathbb T) \cong L^1_{n \varLambda}(\mathbb T)$ for every 
integer $n \neq 0$.

\section{Rich subspaces}
\label{SectionRichSubspaces}

\begin{definition}
    Let $X$ be a Banach space with the Daugavet property and let $E$ be an
    arbitrary Banach space. An operator $T \in L(X,E)$ is called
    \emph{narrow}, if for every two elements $x, y \in S_X$,
    for every $x^* \in X^*$, and for every $\varepsilon > 0$ there is an 
    element $z \in S_X$ such that $\norm{T(y - z)} + |x^*(y - z)| \leq 
    \varepsilon$ and $\norm{x + z} \geq 2 - \varepsilon$. A closed subspace
    $Y$ of $X$ is said to be \emph{rich}, if the quotient map
    $\pi : X \rightarrow X/Y$ is narrow.
\end{definition}
A rich subspace inherits the Daugavet property. But even a little bit more
is true \cite{KadetsShvidkoyWernerNarrowOperators}*{Theorem~5.2}.
\begin{prop}
    \label{PropositionRichSubspace}
    Let $X$ be a Banach space with the Daugavet property and let $Y$ be a 
    rich subspace. Then for every $x \in S_X$, $y^* \in S_{Y^*}$, and 
    $\varepsilon > 0$ there is some $y \in S_Y$ with $\Re y^*(y) \geq 1 - 
    \varepsilon$ and $\norm{x + y} \geq 2 - \varepsilon$.
    \begin{proof}
        Fix $x \in S_X$, $y^* \in S_{Y^*}$, and $\varepsilon > 0$. 
        Choose $\delta > 0$ with $\frac{1 - 3\delta}{1 + \delta}\geq 1 - 
        \varepsilon$ and$z \in S_Y$ with $\Re y^*(z) \geq 1 - \delta$. 
        Since $Y$ is a rich subspace of $X$, there exists $x_0 \in S_X$ with
        $d(x_0,Y) = d(z - x_0,Y) < \delta$, $|y^*(z - x_0)| \leq \delta$ and 
        $\norm{x + x_0} \geq 2 - \delta$. Fix $y_0 \in Y$ with 
        $\norm{x_0 - y_0} \leq \delta$ and set $y = \frac{y_0}{\norm{y_0}}$.
        Then 
        \begin{align*}
            \Re y^*(y_0) &\geq \Re y^*(z) - |y^*(z-x_0)| - \norm{x_0-y_0}
                \geq 1 - 3\delta\\
            \shortintertext{and}
            \norm{x_0 - y} &\leq \norm{x_0 - y_0} + \norm{y_0 - y} 
                \leq 2 \delta.
        \end{align*}
        So we get by our choice of $\delta$ that $\Re y^*(y) \geq 1 -
        \varepsilon$ and $\norm{x + y} \geq 2 - \varepsilon$.
    \end{proof}
\end{prop}
Let us recall the following characterizations of narrow operators on $C(K)$
spaces \cite{KadetsShvidkoyWernerNarrowOperators}*{Theorem~3.7} and on
$L^1(\varOmega,\varSigma,\mu)$ spaces
\citelist{\cite{KadetsShvidkoyWernerNarrowOperators}*{Theorem~6.1}
\cite{KadetsKaltonWernerRemarks}*{Theorem~2.1}}.
\begin{prop}
    \label{PropositionCharakterizationNarrowC}
    Let $K$ be a compact space without isolated points and let $E$ be a
    Banach space. An operator $T \in L(C(K),E)$ is narrow if and only if
    for every non-empty open set $O$ and every $\varepsilon > 0$ there is 
    a function $f \in S_{C(K)}$ with $\restr{f}{K \setminus O} = 0$ and 
    $\norm{T(f)} \leq \varepsilon$.
\end{prop}
\begin{rem}
    In Proposition~\ref{PropositionCharakterizationNarrowC} the function
    $f$ can be chosen to be real-valued and non-negative. This was proven
    for $C(K,\mathbb R)$ in 
    \cite{KadetsPopovDaugavetPropertyNarrowOperatorsRichSubspaces}*{Lemma~1.4}.
    The same proof works with minor modifications for $C(K, \mathbb C)$ as
    well. 
\end{rem}
\begin{prop}
    \label{PropositionCharakterizationNarrowL1}
    Let $(\varOmega, \varSigma, \mu)$ be a non-atomic probability space and 
    let $E$ be a Banach space. A function $f \in L^1(\varOmega)$ is said to be a
    \emph{balanced $\varepsilon$-peak} on $A \in \varSigma$, if $f$ is
    real-valued, $f \geq -1$, $\supp f \subset A$, 
    $\int_\varOmega f \,d\mu = 0$, and $\mu(\{ f = -1 \}) \geq \mu(A) - 
    \varepsilon$. An operator $T \in L(L^1(\varOmega),E)$ is narrow if
    and only if for every $A \in \varSigma$ and every $\delta, \varepsilon > 0$
    there is a balanced $\varepsilon$-peak $f$ on $A$ with $\norm{T(f)} \leq
    \delta$.    
\end{prop}
\begin{cor}
    \label{CorollaryPeakRichSubspaceC}
    If $C_\varLambda(G)$ is a rich subspace of $C(G)$, then there exists for
    every $x \in G$, every open neighborhood $O$ of $e_G$, and every
    $\varepsilon > 0$ a real-valued and non-negative $f \in S_{C(G)}$ with
    $f(x) = 1$, $\restr{f}{G \setminus (x + O)} = 0$, and 
    $d(f, C_\varLambda(G)) \leq \varepsilon$.
    \begin{proof}
        Let $V$ be a symmetric open neighborhood of $e_G$ with $V + V \subset
        O$. Since $C_\varLambda(G)$ is a rich subspace of $C(G)$, we can pick
        a real-valued, non-negative $g \in S_{C(G)}$ with 
        $\restr{g}{G \setminus V} = 0$ and $d(g, C_\varLambda(G)) \leq
        \varepsilon$. Fix $x_0 \in V$ with $g(x_0) = 1$ and set $f =
        g_{x - x_0}$.
        This function is still at a distance of at most $\varepsilon$ from
        $C_\varLambda(G)$ because $C_\varLambda(G)$ is translation-invariant.
        Furthermore, $f(x) = 1$ and 
        $\restr{f}{G \setminus (x + O)} = 0$ by our choice of $V$.
        In fact, if we pick $y \in G$ with $f(y) \neq 0$, we get that
        $g(y - x + x_0) = f(y) \neq 0$. Consequently, $y - x + x_0 \in V$
        and $y \in x - x_0 + V \subset x + V + V \subset x + O$.
    \end{proof}
\end{cor}

We have seen in Proposition~\ref{PropositionRichSubspace} that a rich subspace 
inherits the Daugavet property. But even more is true. A closed subspace $Y$ of 
$X$ is rich if and only if every closed subspace
$Z$ of $X$ with $Y \subset Z \subset X$ has the Daugavet property
\cite{KadetsShvidkoyWernerNarrowOperators}*{Theorem~5.12}. In order to prove
that a translation-invariant subspace $Y$ of $C(G)$ or $L^1(G)$ is rich, we 
do not have to consider all subspaces of $C(G)$ or $L^1(G)$ containing $Y$ but
only the translation-invariant ones.
\begin{lem}
    \label{LemmaApproximationDaugavetProperty}
    Let $X$ be a Banach space. Suppose that for every $\varepsilon > 0$ there
    is a Banach space $Y$ with the Daugavet property and a surjective operator
    $T : X \rightarrow Y$ with $(1 - \varepsilon) \norm{x} \leq \norm{T(x)} 
    \leq (1 + \varepsilon) \norm{x}$ for all $x \in X$. Then $X$ has the
    Daugavet property.
    \begin{proof}
        Let $S : X \rightarrow X$ be an operator of rank one. We have to show
        that $\norm{\mathrm{Id}_X + S} = 1 + \norm{S}$.
        
        Fix $\varepsilon > 0$. By assumption, there exists a Banach space $Y$
        with the Daugavet property and a surjective operator $T : X \rightarrow
        Y$ with $(1 - \varepsilon) \norm{x} \leq \norm{Tx} \leq (1 +
        \varepsilon) \norm{x}$ for all $x \in X$. It is easy to check
        that for every continuous operator $R : X \rightarrow X$ the norm 
        of $TRT^{-1}$ can be estimated by
        \begin{equation*}
            \frac{1 - \varepsilon}{1 + \varepsilon} \norm{R} \leq
            \norm{TRT^{-1}} \leq \frac{1 + \varepsilon}{1 - \varepsilon}
            \norm{R}.
        \end{equation*}
        Using this estimation and the fact that $Y$ has the Daugavet property,
        we get
        \begin{align*}
            \norm{\mathrm{Id}_X + S} &\geq 
                \frac{1 - \varepsilon}{1 + \varepsilon} \norm{\mathrm{Id}_Y + 
                TST^{-1}}\\
            &=\frac{1 - \varepsilon}{1 + \varepsilon} \left( 1 +
                \norm{TST^{-1}} \right)\\
            &\geq \frac{1 - \varepsilon}{1 + \varepsilon} 
                \left( 1 + \frac{1 - \varepsilon}{1 + \varepsilon} \norm{S}
                \right)
        \end{align*}
        This finishes the proof because $\varepsilon > 0$ was arbitrarily
        chosen.
    \end{proof}
\end{lem}
\begin{prop}
    \label{PropositionRichTranslationInvariant}
    Suppose that $\varLambda$ is a subset of $\varGamma$ such that 
    $C_\varTheta(G)$ has the Daugavet property for all $\varLambda \subset 
    \varTheta \subset \varGamma$. Then $C_\varLambda(G)$ is a rich subspace of 
    $C(G)$. The analogous statement is valid for subspaces of $L^1(G)$.
    \begin{proof}
        We will only prove the result for subspaces of $C(G)$. The proof
        for subspaces of $L^1(G)$ works the same way.
        
        It suffices to show that for arbitrary $f_1, f_2 \in S_{C(G)}$ the
        linear span of $C_\varLambda(G)$, $f_1$ and $f_2$ has the Daugavet
        property \cite{KadetsShvidkoyWernerNarrowOperators}*{Lemma~5.6}. In
        order to do this, we are going to prove that $X = \lin \{ 
        C_\varLambda(G) \cup \{ f_1, f_2 \} \}$ meets the assumptions of 
        Lemma~\ref{LemmaApproximationDaugavetProperty}.

        Fix $\varepsilon > 0$ and let us suppose that 
        $f_1 \notin C_\varLambda(G)$ and $f_2 \notin
        \lin \{ C_\varLambda(G) \cup \{ f_1 \} \}$; the other cases can be 
        treated similarly. Then $X$ is isomorphic to $C_\varLambda(G) \oplus_1 
        \lin \{ f_1 \} \oplus_1 \lin \{ f_2 \}$ and there exists $M > 0$ with
        \begin{equation*}
            M (\norm{h}_\infty + |\alpha| + |\beta|) 
            \leq \norm{h + \alpha f_1 + \beta f_2}_\infty
            \quad (h \in C_\varLambda(G), \alpha, \beta \in \mathbb C).
        \end{equation*}
        Since the trigonometric polynomials are dense in $C(G)$, we can choose
        $g_1, g_2 \in S_{T(G)}$ with $\norm{f_k - g_k}_\infty \leq 
        M \varepsilon$ for $k = 1, 2$. If we define $T: X \rightarrow
        \lin \{ C_\varLambda(G) \cup \{ g_1, g_2 \} \}$ by
        \begin{gather*}
            T(h + \alpha f_1 + \beta f_2) = h + \alpha g_1 + \beta g_2
                \quad (h \in C_\varLambda(G), \alpha, \beta \in \mathbb C),
            \shortintertext{then $T$ is surjective and meets the 
                assumption of Lemma~\ref{LemmaApproximationDaugavetProperty} 
                since}
            \norm{T(h + \alpha f_1 + \beta f_2) - (h + \alpha f_1 + 
                \beta f_2)}_\infty \leq M \varepsilon (|\alpha| + |\beta|) \leq
                \varepsilon \norm{h + \alpha f_1 + \beta f_2}_\infty
        \end{gather*}
        for $h \in C_\varLambda(G)$ and $\alpha, \beta \in \mathbb C$.
        
        To complete the proof, we have to show that 
        $Y = \lin \{ C_\varLambda(G) \cup \{ g_1, g_2 \} \}$ has the Daugavet
        property. Set $\varDelta = \spec g_1 \cup \spec g_2$. 
        Since $g_1$ and $g_2$ are trigonometric polynomials,
        the set $\varDelta$ is finite. By assumption, $C_{\varLambda \cup
        \varDelta}(G)$ has the Daugavet property. The space $Y$ is a
        finite-codimensional subspace of $C_{\varLambda \cup \varDelta}(G)$ and
        has therefore the Daugavet property as well 
        \cite{KadetsShvidkoySirotkinWernerDaugavetProperty}*{Theorem~2.14}.
    \end{proof}
\end{prop}
Not all translation-invariant subspaces of $C(G)$ or $L^1(G)$ which have the
Daugavet property must be rich. The subspace $C_{2\mathbb Z}(\mathbb T)$ has 
the Daugavet property because $C(\mathbb T) \cong C_{2\mathbb Z}(\mathbb T)$
by Corollary~\ref{CorollaryHomomorphism}. But every $f \in 
C_{2\mathbb Z}(\mathbb T)$ satisfies
\begin{equation*}
    f(t) = f(-t) \quad (t \in \mathbb T)
\end{equation*}
and therefore $C_{2\mathbb Z}(\mathbb T)$ cannot be a rich subspace of
$C(\mathbb T)$. Similarly, $L^1_{2\mathbb Z}(\mathbb T)$ has the Daugavet 
property but is not a rich subspace of $L^1(\mathbb T)$.

If $X$ is a Banach space with the Daugavet property, then all operators on $X$
which do not fix $\ell^1$ are narrow
\cite{KadetsShvidkoyWernerNarrowOperators}*{Theorem~4.13}. This implies that 
$Y$ is a rich subspace of $X$, if the quotient space $X/Y$ contains no copy of 
$\ell^1$ or if $(X/Y)^*$ has the Radon-Nikod\'ym property
\cite{KadetsShvidkoyWernerNarrowOperators}*{Proposition~5.3}. Let us apply 
these results to translation-invariant subspaces of $C(G)$ or $L^1(G)$.
\begin{definition}
    \mbox{}
    \begin{enumerate}[\upshape (a)]
        \item $\varLambda$ is called a \emph{Rosenthal set}, if every 
            equivalence class of $L^\infty_\varLambda(G)$ contains a continuous
            member, i.e., if $L^\infty_\varLambda(G) = C_\varLambda(G)$.
        \item $\varLambda$ is called a \emph{Riesz set}, if every $\mu \in
            M_\varLambda(G)$ is absolutely continuous with respect to the Haar
            measure, i.e., if $M_\varLambda(G) = L^1_\varLambda(G)$.
    \end{enumerate}
\end{definition}
Every Sidon set is a Rosenthal set. (Recall that $\varLambda \subset \varGamma$ 
is said to be a \emph{Sidon set}, if there exists a constant $M > 0$ such that
$\sum_{\gamma \in \varLambda} |\hat f(\gamma)| \leq M \norm{f}_\infty$ for all
$f \in T_\varLambda(G)$.) But $\bigcup_{n = 1}^\infty (2n)! \{ 1, \dotsc, 2n \}$
is an example of a Rosenthal set which is not a Sidon set
\cite{RosenthalTrigonometricSeriesAssociatedWeak*ClosedSubspaces}*{Corollary~4}.
Every Rosenthal set is a Riesz set
\cite{LustEnsemblesRosenthalRiesz}*{Th\'eor\`eme~3} and it is a classical 
result due to F. and M.~Riesz that $\mathbb N$ is a Riesz set
\cite{RudinRealComplexAnalysis}*{Theorem~17.13}.
\begin{prop}
    \label{PropositionRieszRosenthal}
    If $\varLambda$ is a Riesz set,then 
    $C_{\varGamma \setminus \varLambda^{-1}}(G)$ is a rich subspace of $C(G)$, 
    and if $\varLambda$ is a Rosenthal set, then 
    $L^1_{\varGamma \setminus \varLambda^{-1}}(G)$ is a rich subspace of 
    $L^1(G)$.
    \begin{proof}
        Suppose that $\varLambda$ is a Riesz set. Since 
        $T_{\varGamma \setminus \varLambda^{-1}}(G)$ is dense in $C_{\varGamma 
        \setminus \varLambda^{-1}}(G)$, we have that $C_{\varGamma \setminus 
        \varLambda^{-1}}(G)^\perp = M_\varLambda(G)$. Hence the dual space of 
        $C(G)/C_{\varGamma \setminus \varLambda^{-1}}(G)$ can be identified with
        $M_\varLambda(G)$. Since $\varLambda$ is a Riesz set, 
        $M_\varLambda(G)$ has the Radon-Nikod\'ym property
        \cite{LustEnsemblesRosenthalRiesz}*{Th\'eor\`eme~2} and 
        $C_{\varGamma \setminus \varLambda^{-1}}(G)$ is a rich subspace of 
        $C(G)$ \cite{KadetsPopovDaugavetPropertyNarrowOperatorsRichSubspaces}*{
        Proposition~5.3}.
    
        Suppose now that $\varLambda$ is a Rosenthal set. We apply the same
        reasoning as before and use the fact that $L^\infty_\varLambda(G)$ has
        the Radon-Nikod\'ym property, if $\varLambda$ is a Rosenthal set
        \cite{LustEnsemblesRosenthalRiesz}*{Th\'eor\`eme~1}.        
    \end{proof}
\end{prop}

In Section~\ref{SectionProducts}, we will give an example of a 
non-Rosenthal set $\varLambda$ such that 
$L^1_{\varGamma \setminus \varLambda^{-1}}(G)$ is a rich subspace of $L^1(G)$. 
In the case of translation-invariant subspaces of $C(G)$, the previous result 
can be strengthened.
\begin{definition}
    A measure $\mu \in M(G)$ is said to be \emph{diffuse} or
    \emph{non-atomic}, if $\mu(B) = 0$ for all countable sets $B \subset G$.
    We denote by $M_{\mathrm{diff}}(G)$ the set of all diffuse members of
    $M(G)$. A subset $\varLambda$ of $\varGamma$ is called a 
    \emph{semi-Riesz set}, if every $\mu \in M_\varLambda(G)$ is diffuse.    
\end{definition}
If $G$ is infinite, the Haar measure on $G$ is diffuse and every Riesz set
of $\varGamma$ is a semi-Riesz set. The set $\{ \sum_{k = 0}^n \varepsilon_k
4^k : n \in \mathbb N, \varepsilon_k \in \{ -1, 0, 1 \} \}$ is an example
of a proper semi-Riesz set 
\cite{WernerDaugavetEquationFunctionSpaces}*{p.~126}.

D.~Werner showed that $C_{\varGamma \setminus \varLambda^{-1}}(G)$ has the 
Daugavet property, if $\varLambda$ is a semi-Riesz set
\cite{WernerDaugavetEquationFunctionSpaces}*{Theorem~3.7}. Combining this
result with the fact that every subset of a semi-Riesz set is still a 
semi-Riesz set, we get by Proposition~\ref{PropositionRichTranslationInvariant}
the following corollary.
\begin{cor}
    If $\varLambda$ is a semi-Riesz set,then 
    $C_{\varGamma \setminus \varLambda^{-1}}(G)$ is a rich subspace of $C(G)$.
\end{cor}
The converse implication is also valid.
\begin{theo}
    \label{TheoremRichSubspaceSemiRiesz}
    If $C_{\varLambda}(G)$ is a rich subspace of $C(G)$, then 
    $C_{\varLambda}(G)^\perp \subset M_{\mathrm{diff}}(G)$.
    \begin{proof}
        Let $[\mu]$ denote the corresponding equivalence class of $\mu$ in
        $M(G)/C_{\varLambda}(G)^\perp$. It suffices to show the following:
        For every $x \in G$, every $\alpha \in \mathbb{C}$, and every
        $\mu \in M(G)$ with $\mu(\{x\}) = 0$ we have
        $\norm{[\alpha \delta_x] + [\mu]} = |\alpha| + \norm{[\mu]}$. 
        Indeed, if the preceding statement is true, we get for every $\mu \in
        C_{\varLambda}(G)^\perp$ and every $x \in G$ that
        \begin{equation*}
            0 = \norm{[\mu]} = \norm{[\mu(\{x\})\delta_x] + 
            [\mu - \mu(\{x\})\delta_x]}
            = |\mu(\{x\})| + \norm{[\mu - \mu(\{x\})\delta_x]}.
        \end{equation*}
        Hence $|\mu(\{x\})| = 0$ and $\mu$ is a diffuse measure.
        
        Fix $x \in G$, $\alpha \in \mathbb C \setminus \{0\}$, $\mu \in M(G)$
        with $\mu(\{x\}) = 0$, and $\varepsilon > 0$. Choose $f \in
        S_{C_\varLambda(G)}$ with $\Re \int_G f \,d\mu \geq \norm{[\mu]} -
        \varepsilon$. Since $|\mu|$ is a regular Borel measure and $f$ is a
        continuous function, there is an open neighborhood $O$ of $e_G$ with 
        \mbox{$|\mu|(x + O) < \varepsilon$} and 
        $|f(x) - f(x + y)| < \varepsilon$ for 
        all $y \in O$. As $C_\varLambda(G)$ is a rich subspace of $C(G)$, we can
        pick by Corollary~\ref{CorollaryPeakRichSubspaceC} a real-valued, 
        non-negative $g_0 \in S_{C(G)}$ with \mbox{$g_0(x) = 1$}, 
        $\restr{g_0}{G \setminus (x + O)} = 0$ and $d(g_0, C_\varLambda(G)) < 
        \varepsilon$. Let $g$ be an element of $C_\varLambda(G)$ with 
        $\norm{g - g_0}_\infty \leq \varepsilon$.
        
        If we set $h_0 = f + \left(\frac{|\alpha|}{\alpha} - f(x) \right)
        g_0$ and $h = f + \left(\frac{|\alpha|}{\alpha} - f(x) \right) g$,
        then $h \in C_\varLambda(G)$ and $\norm{h - h_0}_\infty \leq 2
        \varepsilon$. Furthermore, 
        \begin{gather}
            \label{TheoremRichSemiRieszEqu1}
            \alpha h_0(x) = |\alpha| 
            \shortintertext{and}
            \label{TheoremRichSemiRieszEqu2}
            \begin{align}
                \Re \int_G h_0 \, d\mu &= 
                    \Re \int_G \left(f + \left(\frac{|\alpha|}{\alpha} -
                    f(x) \right) g_0 \right) \, d\mu\\
                &\geq \norm{[\mu]} - \varepsilon - 2\int_G g_0 \, d|\mu|
                    \nonumber\\
                &= \norm{[\mu]} - \varepsilon - 2\int_{x + O} g_0 \,
                    d|\mu| \nonumber\\
                &\geq \norm{[\mu]} - \varepsilon -2|\mu|(x + O) \nonumber\\
                &\geq \norm{[\mu]} - 3\varepsilon. \nonumber
            \end{align}
        \end{gather}
        Let us estimate the norm of $g$.
        We get for $y \in G \setminus (x + O)$
        \begin{align*}
            |h(y)| &= \left| f(y) + \left( \frac{|\alpha|}{\alpha} - f(x)
                \right) g(y) \right|
                \leq \norm{f}_\infty + 2 \lVert \restr{g}{G \setminus (x + O)}
                \rVert_\infty \leq 1 + 2 \varepsilon
            \shortintertext{and for $y \in x + O$}
            |h(y)| &= \left| f(y) + \left( \frac{|\alpha|}{\alpha} - f(x) 
                \right) g(y) \right|
                \leq |f(y) - f(x)g_0(y)| + g_0(y) + 2\norm{g - g_0}_\infty\\
            &\leq |f(y) - f(x)| + |f(x)|(1 - g_0(y)) + g_0(y) + 2 \varepsilon \\
            &\leq \varepsilon + (1 - g_0(y)) + g_0(y) + 2 \varepsilon
                =1 + 3\varepsilon.
        \end{align*}
        Hence $\norm{h}_\infty \leq 1 + 3\varepsilon$. Combining this estimate
        with (\ref{TheoremRichSemiRieszEqu1}) and
        (\ref{TheoremRichSemiRieszEqu2}), we get
        \begin{align*}
        (1 + 3\varepsilon)\norm{[\alpha \delta_x] + [\mu]} &\geq
            \left| \int_G h \,d(\alpha \delta_x + \mu)\right|\\
        &\geq \left| \int_G h_0\,d(\alpha \delta_x + \mu)\right|
            - 2\varepsilon \norm{\alpha \delta_x + \mu}\\
        &\geq |\alpha| + \norm{[\mu]} - 3\varepsilon
            - 2\varepsilon \norm{\alpha \delta_x + \mu}.
        \end{align*}
        We can choose $\varepsilon > 0$ arbitrarily small and so 
        $\norm{[\alpha \delta_x] + [\mu]} = |\alpha| + \norm{[\mu]}$.
    \end{proof}
\end{theo}
\begin{cor}
    \label{CorollaryCharacterizationRichC}
    The space $C_\varLambda(G)$ is a rich subspace of $C(G)$ if and only if 
    $\varGamma \setminus \varLambda^{-1}$ is a semi-Riesz set.
\end{cor}

A linear projection $P$ on a Banach space $X$ is called an
\emph{$L$-projection}, if
\begin{equation*}
    \norm{x} = \norm{P(x)} + \norm{x - P(x)} \quad (x \in X).
\end{equation*} 
A closed subspace of $X$ is called an \emph{$L$-summand}, if it is the
range of an $L$-projection. In the proof of 
Theorem~\ref{TheoremRichSubspaceSemiRiesz} we showed that every Dirac measure
$\delta_x$ still has norm one and still spans an $L$-summand, if we 
consider it as an element of $C_\varLambda(G)^*$. Such subspaces are
called \emph{nicely embedded} and were studied by D.~Werner
\cite{WernerDaugavetEquationFunctionSpaces}. His proof of the fact
that $C_\varLambda(G)$ has the Daugavet property, if $\varGamma \setminus
\varLambda^{-1}$ is a semi-Riesz set, is as well based on the observation
that then $C_\varLambda(G)$ is nicely embedded.

Let us present an alternative proof of 
Corollary~\ref{CorollaryCharacterizationRichC} for the case that $G$ is 
metrizable. It is based on results of V.\,M.~Kadets and M.\,M.~Popov 
\cite{KadetsPopovDaugavetPropertyNarrowOperatorsRichSubspaces}.
We say that an operator $T \in L(C(G),E)$ vanishes at a 
point $x \in G$ and write $x \in \van T$, if  there exists a sequence
$\left(O_n\right)_{n \in \mathbb{N}}$ of open neighborhoods of $x$ with
\mbox{$\operatorname{diam}{O_n} \longrightarrow 0$} and a sequence 
$\left(f_n\right)_{n \in \mathbb{N}}$ of non-negative functions satisfying
that \mbox{$f_n \in S_{C(K)}$, $\restr{f_n}{G \setminus O_n} = 0$}, $\left(f_n
\right)_{n \in \mathbb{N}}$ converges pointwise to $\chi_{\{x\}}$, and 
$\norm{T(f_n)} \longrightarrow 0$. An operator $T$ is narrow if and only if
$\van T$ is dense in $G$ 
\cite{KadetsPopovDaugavetPropertyNarrowOperatorsRichSubspaces}*{Lemma~1.6}.
Furthermore, $x \in \van T$ if and only if for any functional $e^* \in E^*$
the point $x$ is not an atom of the measure corresponding to $T^*(e^*)$
\cite{KadetsPopovDaugavetPropertyNarrowOperatorsRichSubspaces}*{Lemma~1.7}.
Let $\varLambda$ be a subset of $\varGamma$, let $\pi: C(G) \rightarrow
C(G)/C_{\varLambda}(G)$ be the canonical quotient map and note that
\begin{equation*}
    \operatorname{ran}(\pi^*) =  
    C_{\varLambda}(G)^\perp = M_{\varGamma \setminus \varLambda^{-1}}(G).
\end{equation*}
If $\varGamma \setminus \varLambda^{-1}$ is a semi-Riesz set, then every element
of $M_{\varGamma \setminus \varLambda^{-1}}(G)$ is a diffuse measure. Therefore
$\van \pi = G$ and $\pi$ is a narrow operator. Conversely, if $\pi$ is
narrow, it is an easy consequence of Corollary~\ref{CorollaryPeakRichSubspaceC}
that $\van \pi = G$. Therefore, $M_{\varGamma \setminus \varLambda^{-1}}(G)$ 
must consist of diffuse measures and $\varGamma \setminus \varLambda^{-1}$ is a 
semi-Riesz set.

Let $\varLambda$ be a subset of $\mathbb Z$ and let $\lambda_1, \lambda_2,
\dotsc$ be an enumeration of $\varLambda$ with $|\lambda_1| \leq |\lambda_2|
\leq \dotsb$. We say that $\varLambda$ is \emph{uniformly distributed}, if
\begin{equation*}
\frac{1}{n} \sum_{k = 1}^n t^{\lambda_k} \longrightarrow 0
\quad(t \in \mathbb T, t \neq 1).
\end{equation*}
R.~Demazeux proved that $C_\varLambda(\mathbb T)$ is a rich subspace of
$C(\mathbb T)$, if $\varLambda$ is uniformly distributed
\cite{DemazeuxCentresDaugavetOperateursComposition}*{Th\'eor\`eme~I.1.7}.
Corollary~\ref{CorollaryCharacterizationRichC} shows that
$\mathbb Z \setminus (-\varLambda)$ is a semi-Riesz set, if $\varLambda$ is
uniformly distributed.
\begin{theo}
    If $L^1_{\varGamma \setminus \varLambda^{-1}}(G)$ is a rich subspace of 
    $L^1(G)$, then $\varLambda$ is a semi-Riesz set.
    \begin{proof}
    The following proof is based on arguments used by G.~Godefroy, N.\,J.~Kalton
    and D.~Li
    \cite{GodefroyKaltonLiOperatorsBetweenSubspacesQuotients}*{
    Proposition~III.10} and N.\,J.~Kalton 
    \cite{KaltonEndomorphismsLp}*{Theorem~5.4}.
    
    Suppose that $\varLambda$ is not a semi-Riesz set. We will
    show that $L^1_{\varGamma \setminus \varLambda^{-1}}(G)$ is not a rich 
    subspace of $L^1(G)$.
    
    Let $\mu \in M_\varLambda(G)$ be a non-diffuse measure.
    We may assume that $\mu = \delta_{e_G} + \nu$ with $\nu(\{e_G\}) = 0$.
    (If $\mu$ is not of this form, fix $x \in G$ with $\mu(\{x\}) \neq 0$ and
    consider the measure $\frac 1{\mu(\{x\})}\mu*\delta_{-x} \in 
    M_\varLambda(G)$.) Let $R,S,T : L^1(G) \rightarrow L^1(G)$ be the 
    convolution operators defined by $R(f) = \mu*f$, $S(f) = \nu*f$ and 
    $T(f) = |\nu|*f$. Note that $R = \mathrm{Id} + S$. Recall that for every 
    $\lambda \in M(G)$ and $f \in L^1(G)$ we have
    \begin{equation*}
        (\lambda*f)(x) = \int_G f(x - y) \,d\lambda(y)
    \end{equation*}
    for $m$-almost all $x \in G$ 
    \cite{HewittRossAbstractHarmonicAnalysisI}*{Theorem~V.20.12}. Therefore
    $\norm{S(\chi_E)}_1 \leq \norm{T(\chi_E)}_1$ for all Borel sets $E$ of $G$.
    
    We will first show that there exists $A \in \mathcal B (G)$ 
    with $m(A) > 0$ such that $\restr{R}{L^1(A)}$ is an isomorphism onto its
    image. (We write $L^1(A)$ for the subspace $\{ f \in L^1(G): 
    \supp f \subset A \}$.) Since $\nu(\{e_G\}) = 0$, we can choose a sequence 
    $(O_n)_{n \in \mathbb N}$ of open neighborhoods of 
    $e_G$ with $|\nu|(O_n) \longrightarrow 0$.
    For each $n \in \mathbb N$, use Lemma~\ref{LemmaCovering} to find a covering
    of $G$ by disjoint Borel sets let $B_{n,1}, \dotsc, B_{n,N_n}$ with 
    $B_{n,k} - B_{n,k} \subset O_n$ for $k = 1, \dotsc, N_n$. 
    Set for every $n \in \mathbb N$
    \begin{equation*}
    R_n = \sum_{k = 1}^{N_n} P_{B_{n,k}} R P_{B_{n,k}},
    \quad
    S_n = \sum_{k = 1}^{N_n} P_{B_{n,k}} S P_{B_{n,k}}, 
    \quad \text{and} \quad
    T_n = \sum_{k = 1}^{N_n} P_{B_{n,k}} T P_{B_{n,k}} 
    \end{equation*}
    where $P_E$ denotes for every $E \in \mathcal B(G)$ the projection 
    from $L^1(G)$ onto $L^1(E)$ defined by $P_E(f) = \chi_E f$. Let for every
    $n \in \mathbb N$ the map $\rho_n$ be defined by
    \begin{equation*}
        \rho_n(E) = \norm{T_n (\chi_E)}_1 \quad (E \in \mathcal B(G)).
    \end{equation*}
    Since $T_n$ is continuous and maps positive functions to positive 
    functions, it is a consequence of the monotone convergence theorem
    that $\rho_n$ is a positive Borel measure on $G$. Every $\rho_n$ is 
    absolutely continuous with respect to the Haar measure $m$ and has 
    Radon-Nikod\'ym derivative $\omega_n$. For each $n \in 
    \mathbb N$, we get
    \begin{align*}
        \rho_n(G) &= \norm{T_n (\chi_G)}_1 = \sum_{k = 1}^{N_n} 
            \int_{B_{n,k}} T (\chi_{B_{n,k}})(x) \, dm(x)\\
        &= \sum_{k = 1}^{N_n} \int_{B_{n,k}} \int_G \chi_{B_{n,k}}(x - y) \,
            d|\nu|(y) dm(x)\\
        &= \sum_{k = 1}^{N_n} \int_{B_{n,k}} |\nu|(x - B_{n,k}) \,dm(x)\\
        &\leq \sum_{k = 1}^{N_n} \int_{B_{n,k}} |\nu|(B_{n,k} - B_{n,k})
        \,dm(x) \leq |\nu|(O_n).
    \end{align*}
    Therefore, $\rho_n(G) \longrightarrow 0$ and in particular $\omega_n
    \longrightarrow 0$ in $m$-measure. So there exists a Borel set $B_0$ of
    $G$ with $m(B_0) > 0$ and $n_0 \in \mathbb N$ satisfying
    \begin{equation*}
        \omega_{n_0}(x) \leq \frac{1}{2} \quad(x \in B_0).
    \end{equation*}
    Consequently, $\norm{S_{n_0}(\chi_E)}_1 \leq \norm{T_{n_0} (\chi_E)}_1 
    \leq \frac 12 m(E)$ for
    all Borel sets $E \subset B_0$ and $\lVert \restr{S_{n_0}}{L^1(B_0)}
    \rVert \leq \frac 12$. Therefore $\restr{(\mathrm{Id} + S_{n_0})}{L^1(B_0)}
    = \restr{R_{n_0}}{L^1(B_0)}$
    is an isomorphism onto its image. Fix $k_0 \in \{1, \dotsc, N_{n_0}\}$ with
    $m(B_0 \cap B_{n_0,k_0}) > 0$ and set \mbox{$A = B_0 \cap B_{n_0,k_0}$}. 
    Then $\restr{R}{L^1(A)}$ is an isomorphism because $\norm{R_{n_0}(f)}_1 \leq
    \norm{R(f)}_1$ for all $f \in L^1(G)$.
    
    We will now finish the proof by showing that $L^1_{\varGamma \setminus
    \varLambda^{-1}}(G)$ is not a rich subspace of $L^1(G)$. Let $\pi: L^1(G)
    \rightarrow L^1(G)/\ker R$ be the canonical quotient map and let 
    $\widetilde R : L^1(G)/\ker R \rightarrow L^1(G)$ be a bounded operator
    with $R = \widetilde R \circ \pi$. Since $\restr{R}{L^1(A)}$ is an
    isomorphism, $\restr{\pi}{L^1(A)}$ is bounded from below. By
    Proposition~\ref{PropositionCharakterizationNarrowL1}, $\pi$ cannot be
    a narrow operator. $L^1_{\varGamma \setminus \varLambda^{-1}}(G)$ is 
    contained in $\ker R$ and is therefore not a rich subspace of $L^1(G)$.
    \end{proof}
\end{theo}
\begin{cor}
    \label{CorollaryL1rich}
    If $L^1_\varLambda(G)$ is a rich subspace of $L^1(G)$, 
    then $C_\varLambda(G)$ is a rich subspace of $C(G)$.
\end{cor}
The space $C_\mathbb N(\mathbb T)$ is a rich subspace of $C(\mathbb T)$, but
$L^1_\mathbb N(\mathbb T)$ has the Radon-Nikod\'ym property and therefore
not the Daugavet property. So the converse of Corollary~\ref{CorollaryL1rich}
is not true.

\section{Products of compact abelian groups}
\label{SectionProducts}

Let $G_1$ and $G_2$ be compact abelian groups with normalized Haar measures 
$m_1$ and $m_2$. The direct product $G = G_1 \times G_2$ is again a compact 
abelian group, if we endow it with the product topology. 
If $f: G_1 \rightarrow \mathbb C$ and $g: G_2 \rightarrow \mathbb C$, we denote 
by $f \otimes g$ the function $(x,y) \mapsto f(x)g(y)$. The dual group of $G$ 
can now be identified with $\varGamma_1 \times \varGamma_2$ because every 
$\gamma \in \varGamma$ is of the form
$\gamma_1 \otimes \gamma_2$ with $\gamma_1 \in \varGamma_1$ and $\gamma_2 \in 
\varGamma_2$ \cite{RudinFourierAnalysis}*{Theorem~2.2.3}. Furthermore, the Haar 
measure on $G$ coincides with the product measure $m_1 \times m_2$ 
\cite{HewittRossAbstractHarmonicAnalysisI}*{Example~IV.15.17.(i)}.
\begin{prop}
    \label{PropositionProductOfRichRichC}
    Let $G_1$ be an infinite, compact abelian group, let $G_2$ be an arbitrary, 
    compact abelian group, let $\varLambda_1$ be a subset of $\varGamma_1$,
    and let $\varLambda_2$ be a subset of $\varGamma_2$. 
    \begin{enumerate}[\upshape(a)]
        \item Suppose that $C_{\varLambda_1}(G_1)$ is a rich subspace of 
            $C(G_1)$ and that $C_{\varLambda_2}(G_2)$ is a rich subspace of 
            $C(G_2)$ (or, if $G_2$ is finite, that $\varLambda_2 = 
            \varGamma_2$). Then $C_{\varLambda_1 \times \varLambda_2}
            (G_1 \times G_2)$ is a rich subspace of $C(G_1 \times G_2)$.
        \item Suppose that $C_{\varLambda_1}(G_1)$ is a rich subspace of 
            $C(G_1)$ and that $\varLambda_2$ is non-empty. Then
            $C_{\varLambda_1 \times \varLambda_2}(G_1 \times G_2)$ has the 
            Daugavet property.
    \end{enumerate}    
    \begin{proof}
        Set $G = G_1 \times G_2$ and $\varLambda = 
        \varLambda_1 \times \varLambda_2$.
        
        We start with part (a). Let $O$ be a non-empty open set of $G$ and  
        $\varepsilon > 0$. By 
        Proposition~\ref{PropositionCharakterizationNarrowC}, we have to find 
        $f \in S_{C(G)}$ with $\restr{f}{G \setminus O} = 0$ and 
        $d(f, C_{\varLambda}(G)) \leq \varepsilon$.
        
        Pick non-empty open sets $O_1 \subset G_1$ and $O_2 \subset G_2$ with
        $O_1 \times O_2 \subset O$ and $\delta > 0$ with $2\delta + \delta^2 
        \leq \varepsilon$. By assumption, there exist $f_k \in S_{C(G_k)}$ and 
        $g_k \in T_{\varLambda_k}(G_k)$ with $\restr{f_k}{G_k \setminus O_k} = 
        0$ and $\norm{f_k - g_k}_\infty \leq \delta$ for $k = 1,2$. 
        If we set $f = f_1 \otimes f_2$ and $g = g_1 \otimes g_2$, 
        then $f \in S_{C(G)}$, $g \in T_{\varLambda}(G)$, 
        and $\restr{f}{G \setminus O} = 0$. Furthermore,
        \begin{align*}
            d(f, C_{\varLambda}(G)) &\leq
                \norm{f - g}_\infty\\
            &\leq \norm{f_1}_\infty \norm{f_2 - g_2}_\infty +
                \norm{g_2}_\infty \norm{f_1 - g_1}_\infty\\
            &\leq \delta + (1 + \delta) \delta \leq \varepsilon.
        \end{align*}
        
        Let us now consider part (b). The space 
        $C_{\varGamma_1 \times \varLambda_2}(G)$ can canonically be identified 
        with $C(G_1, C_{\varLambda_2}(G_2))$, the space of all continuous 
        functions from $G_1$ into $C_{\varLambda_2}(G_2)$, and has therefore the 
        Daugavet property \cite{KadetsRemarksDaugavetEquation}*{Theorem~4.4}.
        We will prove that $C_{\varLambda}(G)$
        is a rich subspace of $C_{\varGamma_1 \times \varLambda_2}(G)$.
        For this, it is sufficient to show that for every non-empty open 
        set $O$ of $G_1$, every $g \in T_{\varLambda_2}(G_2)$ with 
        $\norm{g}_\infty = 1$, and every
        $\varepsilon > 0$ there exists $f \in S_{C(G_1)}$ with
        $\restr{f}{G_1 \setminus O} = 0$ and $d(f \otimes g, C_{\varLambda}(G)) 
        \leq \varepsilon$
        \cite{BilikKadetsShvidkoySirotkinWernerNarrowOperatorsSupNormedSpaces}*
        {Proposition~4.3.(a)}. Since $C_{\varLambda_1}(G_1)$ is a rich subspace 
        of $C(G_1)$, there exist \mbox{$f \in S_{C(G_1)}$} and $h \in T_{\varLambda_1}
        (G_1)$ with $\restr{f}{G_1 \setminus O} = 0$ and $\norm{f - h}_\infty
        \leq \varepsilon$. Then \mbox{$h \otimes g \in T_{\varLambda}(G)$} and
        \begin{equation*}
            d(f \otimes g, C_{\varLambda}(G))
            \leq \norm{f \otimes g - h \otimes g}_\infty \leq
            \norm{f - h}_\infty \norm{g}_\infty \leq \varepsilon. \qedhere
        \end{equation*}
    \end{proof}
\end{prop}
\begin{prop}
    \label{PropositionFactorRichC}
    Let $G$ be the product of two compact abelian groups $G_1$ and $G_2$ and 
    denote by $p$ the projection from $\varGamma = \varGamma_1 \times 
    \varGamma_2$ onto $\varGamma_1$. If $C_\varLambda(G)$ is a
    rich subspace of $C(G)$, then $C_{p[\varLambda]}(G_1)$ is a rich subspace 
    of $C(G_1)$ (or $p[\varLambda] = \varGamma_1$, if $G_1$ is finite).
    \begin{proof}
        Let $O$ be a non-empty open set of $G_1$ and $\varepsilon > 0$. 
        By Proposition~\ref{PropositionCharakterizationNarrowC}, we have to 
        find $f \in S_{C(G_1)}$ with $\restr{f}{G_1 \setminus O} = 0$ and 
        $d(f, C_{p[\varLambda]}(G_1)) \leq \varepsilon$. (Note that this is 
        sufficient in the case of finite $G_1$ as well.)
        
        Since $C_\varLambda(G)$ is a rich subspace of $C(G)$, there exist
        $f_0 \in S_{C(G)}$ and $g_0 \in T_\varLambda(G)$ with
        $\restr{f_0}{G \setminus (O \times G_2)} = 0$ and 
        $\norm{f_0 - g_0}_\infty \leq \varepsilon$. Fix $(x_0, y_0) \in G$ with 
        $|f_0(x_0,y_0)| = 1$.
        Setting $f = f_0(\,\cdot\,,y_0)$ and $g = g_0(\,\cdot\,,y_0)$, we get
        that $f \in S_{C(G_1)}$, $g \in T_{p[\varLambda]}(G_1)$, and 
        $\restr{f}{G_1 \setminus O} = 0$. Finally,
        \begin{equation*}
            d(f, C_{p[\varLambda]}(G_1)) \leq \norm{f - g}_\infty \leq
            \norm{f_0 - g_0}_\infty \leq \varepsilon. \qedhere
        \end{equation*}
    \end{proof}
\end{prop}
\begin{prop}
    \label{PropositionProductOfRichRichL1}
    Let $G_1$ and $G_2$ be infinite, compact abelian groups, let $\varLambda_1$
    be a subset of $\varGamma_1$, and let $\varLambda_2$ be a subset of 
    $\varGamma_2$. 
    \begin{enumerate}[\upshape(a)]
        \item Suppose that $L^1_{\varLambda_2}(G_2)$ is a rich subspace of 
            $L^1(G_2)$. Then $L^1_{\varGamma_1 \times \varLambda_2}
            (G_1 \times G_2)$ is a rich subspace of $L^1(G_1 \times G_2)$.
        \item Suppose that $L^1_{\varLambda_1}(G_1)$ is a rich subspace of
            $L^1(G_1)$ and that $\varLambda_2$ is non-empty. Then
            $L^1_{\varLambda_1 \times \varLambda_2}(G_1 \times G_2)$ has the 
            Daugavet property.
    \end{enumerate}
    \begin{proof}
        Set $G = G_1 \times G_2$ and 
        $\varLambda = \varLambda_1 \times \varLambda_2$.
        
        We start with part (a). The space $L^1(G)$  can canonically
        be identified with the Bochner space $L^1(G_1, L^1(G_2))$ and
        $L^1_{\varGamma_1 \times \varLambda_2}(G)$ with the subspace
         $L^1(G_1, L^1_{\varLambda_2}(G_2))$. Since $L^1_{\varLambda_2}(G_2)$
        is a rich subspace of $L^1(G_2)$, the space $L^1(G_1, 
        L^1_{\varLambda_2}(G_2))$ is rich in $L^1(G_1, L^1(G_2))$
        \cite{KadetsKaltonWernerRemarks}*{Lemma~2.8}.
        
        Let us now consider part (b). Identifying again $L^1_{\varGamma_1 \times 
        \varLambda_2}(G)$ with the Bochner space $L^1(G_1,
        L^1_{\varLambda_2}(G_2))$, we see that 
        $L^1_{\varGamma_1 \times \varLambda_2}(G)$ has the Daugavet property 
        \cite{KadetsShvidkoySirotkinWernerDaugavetProperty}*{Example after 
        Theorem~2.3}. We will show that $L^1_{\varLambda}(G)$ is a rich subspace
        of $L^1_{\varGamma_1 \times \varLambda_2}(G)$. For this, it is 
        sufficient to find for every Borel set $A$ of $G_1$, every $g \in 
        T_{\varLambda_2}(G_2)$ with $\norm{g}_1 = 1$, and every $\delta, 
        \varepsilon > 0$ a balanced $\varepsilon$-peak $f$ on $A$ with
        $d(f \otimes g, L^1_{\varLambda}(G)) \leq \delta$ 
        \cite{BoykoKadetsWernerNarrowOperatorsBochnerL1}*{Theorem~2.4}. 
        Since $L^1_{\varLambda_1}(G_1)$ is a rich subspace of $L^1(G_1)$,
        there exist a balanced $\varepsilon$-peak $f$ on $A$ and $h \in 
        T_{\varLambda_1}(G_1)$ with $\norm{f - h}_1 \leq \delta$. Then
        $h \otimes g \in T_{\varLambda}(G)$ and
        \begin{equation*}
            d(f \otimes g, L^1_{\varLambda}(G)
            \leq \norm{f \otimes g - h \otimes g}_1 = \norm{f - h}_1 \norm{g}_1
            \leq \delta. \qedhere
        \end{equation*}
    \end{proof}
\end{prop}
\begin{prop}
    \label{PropositionFactorRichL1}
    Let $G$ be the product of two compact abelian groups $G_1$ and $G_2$ and 
    denote by $p$ the projection from $\varGamma = \varGamma_1 \times 
    \varGamma_2$ onto $\varGamma_1$. If $L^1_\varLambda(G)$ is a
    rich subspace of $L(G)$, then $L^1_{p[\varLambda]}(G_1)$ is a rich subspace 
    of $L^1(G_1)$ (or $p[\varLambda] = \varGamma_1$, if $G_1$ is finite).
    \begin{proof}
        If $p[\varLambda] = \varGamma_1$, we have nothing to show. So let us 
        assume that there exists $\gamma \in \varGamma_1 \setminus 
        p[\varLambda]$. The map $f \mapsto (\overline \gamma \otimes 
        \mathbf 1_{G_2}) f$ is an 
        isometry from $L^1(G)$ onto $L^1(G)$ and maps $L^1_\varLambda(G)$ onto
        $L^1_{(\overline \gamma, \mathbf 1_{G_2})\varLambda}(G)$. Analogously,
        the map $f \mapsto \overline \gamma f$ is an isometry from
        $L^1(G_1)$ onto $L^1(G_1)$ and maps $L^1_{p[\varLambda]}(G_1)$ onto
        $L^1_{\overline \gamma p[\varLambda]}(G_1)$. Note that $\overline \gamma
        p[\varLambda] = p[(\overline \gamma, \mathbf 1_{G_2})\varLambda]$ and 
        that $\mathbf 1_{G_1} \notin \overline \gamma p[\varLambda]$. Since
        $L^1_\varLambda(G)$ is a rich subspace of $L^1(G)$ if and only if
        $L^1_{(\overline \gamma, \mathbf 1_{G_2})\varLambda}(G)$ is a rich
        subspace of $L^1(G)$, and since $L^1_{p[\varLambda]}(G_1)$ is a rich 
        subspace of $L^1(G_1)$ if and only if $L^1_{\overline \gamma 
        p[\varLambda]}(G_1)$ is a rich subspace of $L^1(G_1)$, we may assume 
        without loss of generality that $\mathbf 1_{G_1} \notin p[\varLambda]$.
        
        Fix a Borel subset $A$ of $G_1$ and $\delta, \varepsilon > 0$.
        By Proposition~\ref{PropositionCharakterizationNarrowL1}, we have
        to find a balanced $\varepsilon$-peak $f$ on $A$ with 
        $d(f, L^1_{p[\varLambda]}(G_1)) \leq \delta$. By assumption, 
        $L^1_\varLambda(G)$ is a rich subspace of $L^1(G)$ and therefore there 
        exist a balanced $\frac\varepsilon 3$-peak $f_0$ on $A \times G_2$ and
        $g \in T_\varLambda(G)$ with $\norm{f_0 - g}_1 \leq \frac\delta 6$. Set
        \begin{align*}
            B &= \{ y \in G_2 : m_1(\{f_0(\,\cdot\,,y) = -1\}) > m_1(A)
                -\varepsilon \}
            \shortintertext{and}
            C &= \left\{ y \in G_2 : \norm{f_0(\,\cdot\,,y) - g(\,\cdot\,,y)}_1 
                \leq \frac \delta  2 \right\}.
        \end{align*}
        Note that we may assume that $B$ and $C$ are measurable
        \cite{HewittRossAbstractHarmonicAnalysisI}*{Theorem~III.13.8}.
        We then get
        \begin{align*}
            m_1(A) - \frac \varepsilon 3 &\leq m(\{f_0 = -1 \}) = 
                \int_{G_2}\int_{G_1} \chi_{\{f_0 = -1\}}(x,y)\,dm_1(x)dm_2(y)\\
            &= \int_{G_2} m_1(\{f_0(\,\cdot\,,y)= -1\})\,dm_2(y)\\
            &\leq m_2(B)m_1(A) + (1 - m_2(B))(m_1(A) - \varepsilon)\\
            &= m_1(A) + m_2(B)\varepsilon - \varepsilon
            \shortintertext{and}
            \frac \delta 6 &\geq \norm{f_0 - g}_1 =
                \int_{G_2}\norm{f_0(\,\cdot\,,y) - g(\,\cdot\,,y)}_1\,dm_2(y)\\
            &\geq \frac \delta 2(1 - m_2(C)).
        \end{align*}
        Hence $m_2(B) \geq \frac 2 3$ and $m_2(C) \geq \frac 2 3$. Therefore 
        $B \cap C \neq \emptyset$ and we can choose $y_0 \in B \cap C$.
        
        Let us gather the properties of $f_0(\,\cdot\,,y_0) \in L^1(G_1)$. 
        It is clear that $f_0(\,\cdot\,,y_0)$ is real-valued, 
        $f_0(\,\cdot\,,y_0) \geq -1$, and $\supp f_0(\,\cdot\,,y_0) 
        \subset A$. As $y_0 \in B \cap C$, we have  
        $m_1(\{f_0(\,\cdot\,,y_0) = -1\}) > m_1(A) - \varepsilon$ and
         $\norm{f_0(\,\cdot\,,y_0) - g(\,\cdot\,,y_0)}_1 \leq \frac\delta 2$.
        The function $g(\,\cdot\,,y_0)$ belongs to $T_{p[\varLambda]}(G)$ and 
        $\mathbf 1_{G_1} \notin p[\varLambda]$. So $\int_{G_1}
        g(x,y_0)\,dm_1(x) = 0$ and $|\int_{G_1} f_0(x,y_0) \,dm_1(x)| \leq
        \frac\delta 2$. Modifying $f_0(\,\cdot\,,y_0)$ a little bit, we get a
        balanced $\varepsilon$-peak $f$ on $A$ with
        $\norm{f - g(\,\cdot\,,y_0)}_1 \leq \delta$. 
    \end{proof}
\end{prop}
Set $\varLambda = \mathbb Z \times \{0\}$. Then $\varLambda$ is not a Rosenthal 
set because $C_\varLambda(\mathbb T^2) \cong C(\mathbb T)$ contains a copy of 
$c_0$ \cite{LustPiquardBohrLocalProperties}*{Proof of Theorem~3}. 
But $\mathbb Z^2
\setminus (-\varLambda) = \mathbb Z \times (\mathbb Z \setminus \{0\})$ and
$L^1_{\mathbb Z \times (\mathbb Z \setminus \{0\})}(\mathbb T^2)$ is a rich 
subspace of $L^1(\mathbb T^2)$. So the converse of 
Proposition~\ref{PropositionRieszRosenthal} is not true.

Let us come back to examples of translation-invariant subspaces
that have the Daugavet property but are not rich. The examples mentioned 
in Section~\ref{SectionRichSubspaces} are of the following type: We take
a one-to-one homomorphism $H : \varGamma \rightarrow \varGamma$ that is not 
onto. Then $C_{H[\varGamma]}(G)$ and $L^1_{H[\varGamma]}(G)$ have the Daugavet 
property but are not rich subspaces of $C(G)$ or $L^1(G)$. In this case
$\bigcap_{\gamma \in H[\varGamma]} \ker(\gamma)$ contains $\ker(H^*) \neq 
\{e_G\}$. Set $\varLambda = \mathbb Z \times \{1\}$. Using 
Proposition~\ref{PropositionProductOfRichRichC}.(b) and 
\ref{PropositionProductOfRichRichL1}.(b), we see that
$C_\varLambda(\mathbb T^2)$ and $L^1_\varLambda(\mathbb T^2)$ have the Daugavet 
property. But they are not rich subspaces of $C(\mathbb T^2)$ or 
$L^1(\mathbb T^2)$ by Proposition~\ref{PropositionFactorRichC} and
\ref{PropositionFactorRichL1}. Furthermore, 
$\bigcap_{\gamma \in \varLambda} \ker(\gamma) = \{(1,1)\}$.

\section{Quotients with respect to translation-invariant subspaces}
\label{SectionQuotients}

We are going to study  quotients of the form 
$C(G)/C_\varLambda(G)$ and $L^1(G)/L ^1_\varLambda(G)$. The following lemma is
the key ingredient for all results of this section.
\begin{lem}
    \label{LemmaNormRestrictionf}
    If we interpret $f \in C(G)$ as a functional on $M(G)$, we have
    \begin{equation*}
        \lVert\restr{f}{L^1_\varLambda(G)}\rVert = 
        \lVert\restr{f}{M_\varLambda(G)}\rVert.
    \end{equation*}
    Analogously, if we interpret $g \in L^1(G)$ as a functional on 
    $L^\infty(G)$, we have
    \begin{equation*}
        \lVert\restr{g}{C_\varLambda(G)}\rVert = 
        \lVert\restr{g}{L^\infty_\varLambda(G)}\rVert.
    \end{equation*}
    \begin{proof}
        We will just show the first statement. The proof of the second
        statement works the same way.
        
        It is clear that $\lVert\restr{f}{L^1_\varLambda(G)}\rVert \leq
        \lVert\restr{f}{M_\varLambda(G)}\rVert$ because $L^1_\varLambda(G) 
        \subset M_\varLambda(G)$. In order to prove the reverse inequality, we 
        may assume without loss of generality that 
        $\lVert\restr{f}{M_\varLambda(G)}\rVert = 1$. Fix $\varepsilon > 0$ and
        an approximate unit $(v_j)_{j \in J}$ of $L^1(G)$ that
        fulfills the properties listed in 
        Proposition~\ref{PropositionApproximateUnit}. Pick $\mu
        \in M_\varLambda(G)$ with $\norm{\mu} = 1$ and \mbox{$|\int_G f\,d\mu|
        \geq 1 - \frac\varepsilon 2$}. 
        Using that $\hat v_j(\gamma) \longrightarrow 1$ for every $\gamma 
        \in \varGamma$, we can deduce that 
        \begin{equation*}
        \int_G g\,d(\mu*v_j) \longrightarrow \int_G g \, d\mu
        \end{equation*}
        for every $g \in T(G)$. So $\mu$ is the weak*-limit of 
        $(\mu*v_j)_{j \in J}$ because $T(G)$ is dense in $C(G)$. Fix $j_0 \in J$
        with $|\int_G f \,d(\mu*v_{j_0})| \geq 1 - \varepsilon$. Since
        $\mu*v_{j_0} \in L^1_\varLambda(G)$ and $\norm{\mu*v_{j_0}}_1 
        \leq 1$, we have that $\lVert\restr{f}{L^1_\varLambda(G)}\rVert \geq 1
        - \varepsilon$. As $\varepsilon > 0$ was arbitrarily chosen, this 
        finishes the proof.
    \end{proof}
\end{lem}
\begin{theo}
    \label{TheoremQuotientsC}
    If $L^1_{\varGamma \setminus \varLambda^{-1}}(G)$ is a rich subspace of
    $L^1(G)$, then $C(G)/C_\varLambda(G)$ has the Daugavet property.
    \begin{proof}
        Note that $C_\varLambda(G)^\perp = 
        M_{\varGamma \setminus \varLambda^{-1}}(G)$ because $T_\varLambda(G)$ is
        dense in $C_\varLambda(G)$. We can therefore identify the dual space of 
        $C(G)/C_\varLambda(G)$ with 
        $M_{\varGamma \setminus \varLambda^{-1}}(G)$.
        
        Fix $[f] \in C(G)/C_\varLambda(G)$ with $\norm{[f]} = 1$, $\mu \in
        M_{\varGamma \setminus \varLambda^{-1}}(G)$ with $\norm{\mu} = 1$, and
        $\varepsilon > 0$. By Lemma~\ref{LemmaCharacterizationSlices}, we have
        to find $\nu \in M_{\varGamma \setminus \varLambda^{-1}}(G)$ with
        $\norm{\nu} = 1$, $\Re \int_G f\,d\nu \geq 1 - \varepsilon$, and
        $\norm{\mu + \nu} \geq 2 - \varepsilon$. Let $\mu = \mu_s + g\,dm$ be
        the Lebesgue decomposition of $\mu$ where $\mu_s$ and $m$ are singular
        and $g \in L^1(G)$. 
        
        If we interpret $f$ as a functional on $M(G)$, we have by 
        Lemma~\ref{LemmaNormRestrictionf} that
        \begin{equation*}
            \lVert \restr{f}{L^1_{\varGamma \setminus \varLambda^{-1}}(G)} 
            \rVert = \lVert \restr{f}{M_{\varGamma \setminus\varLambda^{-1}}(G)}
            \rVert =1.
        \end{equation*}
        $L^1_{\varGamma \setminus \varLambda^{-1}}(G)$ is a rich subspace of
        $L^1(G)$ and so there exists by 
        Proposition~\ref{PropositionRichSubspace}
        a function $h \in L^1_{\varGamma \setminus \varLambda^{-1}}(G)$ with
        $\norm{h}_1 = 1$, $\Re \int_G f h \,dm \geq 1 - \varepsilon$,
        and \mbox{$\lVert \frac g {\norm{g}_1} + h\rVert_1 \geq 2 - 
        \varepsilon$}. Setting $\nu = h\,dm$, we therefore get 
        \begin{align*}
            \norm{\mu + \nu} &= \norm{\mu_s} + \norm{g + h}_1 =
                \norm{\mu_s} + \norm{\frac{g}{\norm{g}_1} + h - 
                (1 - \norm{g}_1)\frac{g}{\norm{g}_1}}_1\\
            &\geq \norm{\mu_s} + \norm{\frac{g}{\norm{g}_1} + h}_1 -
                (1 - \norm{g}_1)\\
            &\geq \norm{\mu_s} + (2 - \varepsilon) - (1 - \norm{g}_1)\\
            &= \norm{\mu} + 1 - \varepsilon = 2 - \varepsilon. \qedhere
        \end{align*}
    \end{proof}
\end{theo}
\begin{cor}
    If $\varLambda$ is a Rosenthal set, then $C(G)/C_\varLambda(G)$ has the 
    Daugavet property.
\end{cor}
\begin{theo}
    If $C_{\varGamma \setminus \varLambda^{-1}}(G)$ is a rich subspace of 
    $C(G)$, then $L^1(G)/L^1_\varLambda(G)$ has the Daugavet property.
    \begin{proof}
        Let us  begin as in the proof of Theorem~\ref{TheoremQuotientsC}. We
        can identify the dual space of $L^1(G)/L^1_\varLambda(G)$ with
        $L^\infty_{\varGamma \setminus \varLambda^{-1}}(G)$, because 
        $T_\varLambda(G)$ is dense in $L^1_\varLambda(G)$ and therefore 
        $L^1_\varLambda(G)^\perp = L^\infty_{\varGamma \setminus
        \varLambda^{-1}}(G)$.
        
        Fix $[f] \in L^1(G)/L^1_\varLambda(G)$ with $\norm{[f]} = 1$, $g \in
        L^\infty_{\varGamma \setminus \varLambda^{-1}}(G)$ with 
        $\norm{g}_\infty = 1$, and \mbox{$\varepsilon > 0$}. 
        By Lemma~\ref{LemmaCharacterizationSlices}, 
        we have to find 
        $h \in L^\infty_{\varGamma \setminus \varLambda^{-1}}(G)$
        with $\norm{h}_\infty = 1$, $\Re \int_G fh \,dm \geq 1 -
        \varepsilon$, and $\norm{g + h}_\infty \geq 2 - \varepsilon$.
        
        Choose $\delta \in (0,1)$ with $\frac{1 - 5\norm{f}_1 \delta}
        {1 + 3\delta} \geq 1 - \frac\varepsilon 2$, $\eta > 0$ such that 
        $\int_A |f|\,dm \leq \delta$ for all $A \in \mathcal B(G)$ with 
        $m(A) \leq \eta$, and $t \in \mathbb T$ with
        \begin{equation*}
            m\left(\left\{ \Re t^{-1} g \geq 1 - \frac \varepsilon 2  
            \right\}\right) > 0.
        \end{equation*}        
        If we interpret $f$ as a functional on $L^\infty(G)$, we have by 
        Lemma~\ref{LemmaNormRestrictionf} that
        \begin{equation*}
            \lVert \restr{f}{C_{\varGamma \setminus \varLambda^{-1}}(G)} \rVert
            = \lVert \restr{f}{L^\infty_{\varGamma \setminus\varLambda^{-1}}(G)}
            \rVert = 1.
        \end{equation*}
        Pick $h_0 \in C_{\varGamma \setminus \varLambda^{-1}}(G)$ with
        $\norm{h_0}_\infty = 1$ and $\Re \int_G fh_0\,dm \geq 1
        - \delta$. Since $h_0$ is uniformly continuous, there exists
        an open neighborhood $O$ of $e_G$ with
        \begin{equation*}
            |h_0(x) - h_0(y)| \leq \delta \quad (x - y \in O)
        \end{equation*}
        and $m(O) \leq \eta$. By assumption, 
        $C_{\varGamma \setminus \varLambda^{-1}}(G)$ is a rich subspace of 
        $C(G)$ and so there exist by Corollary~\ref{CorollaryPeakRichSubspaceC} 
        a real-valued,non-negative $p_0 \in S_{C(G)}$ with 
        $\restr{p_0}{G \setminus O} = 0$ and $p_0(e_G) = 1$ and
        $p \in C_{\varGamma \setminus \varLambda^{-1}}(G)$ with 
        $\norm{p_0 - p}_\infty \leq \delta$. Then $V = \{ p_0 > 1 - \delta \}$ 
        is an open neighborhood of $e_G$ and $V \subset O$. An easy compactness 
        rgument shows that there exists $x_0 \in G$ with
        \begin{equation*}
            m\left(\left\{ x \in x_0 + V : \Re t^{-1} g(x) \geq 1 - \frac
            \varepsilon 2 \right\}\right) > 0.
        \end{equation*}
        If we set
        \begin{equation*}
            h_1 = h_0 + (t - h_0(x_0))p_{x_0} 
            \quad \text{and} \quad
            h = \frac{h_1}{\norm{h_1}_\infty},            
        \end{equation*}
        then $h$ is normalized and belongs by construction to 
        $C_{\varGamma \setminus \varLambda^{-1}}(G)$. Let us estimate the norm
        of $h_1$. We get for $x \in G \setminus (x_0 + O)$
        \begin{align*}
            |h_1(x)| &=|h_0(x) + (t - h_0(x_0))p(x - x_0)| \leq 
                \norm{h_0}_\infty + 2 
                \lVert \restr{p}{G \setminus O}\rVert_\infty \leq 1 + 2 \delta
            \shortintertext{and for $x \in x_0 + O$}
            |h_1(x)| &= |h_0(x) + (t - h_0(x_0))p(x - x_0)| \\
            &\leq |h_0(x) - h_0(x_0)p_0(x - x_0)| + p_0(x-x_0) + 2 
                \norm{p - p_0}_\infty\\
            &\leq |h_0(x) - h_0(x_0)| + |h_0(x_0)|(1 - p_0(x - x_0)) + 
                p_0(x - x_0) + 2 \delta\\
            &\leq \delta + (1 - p_0(x - x_0)) + p_0(x - x_0) + 2 \delta\\
            &= 1 + 3 \delta.
        \end{align*}
        Consequently, $\norm{h_1}_\infty \leq 1 + 3 \delta$. Let us
        check that $h$ is as desired. We first observe that
        \begin{align*}
            \Re \int_G f h_1 \,dm 
            &\geq \Re \int_G f h_0 \,dm - 2\int_G |f p_{x_0}|\,dm\\
            &\geq (1 - \delta) - 2\int_{x_0 + O} |f|\,dm - 2\norm{f}_1
                \norm{p_0 - p}_\infty \\
            &\geq (1-\delta) - 2\delta - 2\norm{f}_1 \delta = 1 - (3 + 2
                \norm{f}_1) \delta.
        \end{align*}
        Therefore, $\Re \int_G fh \,dm \geq 1 - \varepsilon$ by our
        choice of $\delta$. If $x \in x_0 + V$, we get 
        \begin{align*}
            \Re t^{-1}h_1(x) &\geq \Re t^{-1}h_0(x) + \Re (1 - t^{-1} h_0(x_0))
                p_0(x - x_0) - 2\norm{p_0 - p}_\infty\\
            &\geq \Re t^{-1}h_0(x) + \Re (1 - t^{-1} h_0(x_0))(1 - \delta)
                - 2\delta\\
            &\geq 1 - 3\delta -|h_0(x) - h_0(x_0)| \geq 1 - 4\delta
        \end{align*}
        and hence $\Re t^{-1}h(x) \geq 1 - \frac \varepsilon 2$ by our choice of
        $\delta$. Thus
        \begin{align*}
            m(\{ \Re t^{-1}(g + h) \geq 2 - \varepsilon\}) &\geq 
                m\left(\left\{ \Re t^{-1}g \geq 1 -  \frac \varepsilon 2 
                \right\} \cap \left\{ \Re t^{-1}h \geq 1 - \frac \varepsilon 2 
                \right\}\right)\\
            &\geq m\left(\left\{ \Re t^{-1}g \geq 1 -  \frac \varepsilon 2 
                \right\} \cap (x_0 + V)\right) > 0
        \end{align*}
        and $\norm{g + h}_\infty \geq 2 - \varepsilon$.
    \end{proof}
\end{theo}
\begin{cor}
    If $\varLambda$ is a semi-Riesz set, then $L^1(G)/L^1_\varLambda(G)$ has the 
    Daugavet property.
\end{cor}

\section{Poor subspaces of $L^1(G)$}

In Section~\ref{SectionQuotients}, we have seen some cases in which the
quotient space $L^1(G)/L^1_\varLambda(G)$ has the Daugavet property. 
Recall that a closed subspace $Y$ of a Banach space $X$ with the Daugavet 
property is rich if and only if every closed subspace $Z$ of $X$ with 
$Y \subset Z \subset X$ has the Daugavet property. A similar 
notion for quotients of $X$ was introduced by V.\,M.~Kadets, 
V.~Shepelska, and D.~Werner 
\cite{KadetsShepelskaWernerQuotientsDaugavetProperty}.
\begin{definition}
    Let $X$ be a Banach space with the Daugavet property. A closed subspace
    $Y$ of $X$ is called \emph{poor}, if $X/Z$ has the Daugavet property for
    every closed subspace $Z \subset Y$.
\end{definition}
The poor subspaces of a Banach space with the Daugavet property can be 
described using a generalized concept of narrow operators
\cite{KadetsShepelskaWernerQuotientsDaugavetProperty}. This leads in the
case of $L^1(\varOmega,\varSigma,\mu)$ to the following characterization
\cite{KadetsShepelskaWernerQuotientsDaugavetProperty}*{Corollary~6.6}.
\begin{prop}
    \label{PropositionPoorL1}
    Let $(\varOmega, \varSigma,\mu)$ be a non-atomic probability space. A 
    subspace $X$ of $L^1(\varOmega)$ is poor if and only if for every 
    $A \in \varSigma$ of positive measure and every $\varepsilon > 0$ there 
    exists $f \in S_{L^\infty(\varOmega)}$ with $\supp f \subset A$ and 
    $\norm{\restr{f}{X}} \leq \varepsilon$ where we interpret $f$ as a 
    functional on $L^1(\varOmega)$.
\end{prop}
Using this characterization, we can build a link to a property that was
studied by G.~Godefroy, N.\,J.~Kalton, and D.~Li 
\cite{GodefroyKaltonLiOperatorsBetweenSubspacesQuotients}. In the sequel,
$(\varOmega, \varSigma, \mu)$ denotes a non-atomic probability space and
$P$ the natural projection from $L^1(\varOmega)^{**}$ onto $L^1(\varOmega)$.
For $A \in \varSigma$, we write $L^1(A)$ for the subspace
$\{ f \in L^1(\varOmega): \supp f \subset A \}$ and
$P_A$ for the projection from $L^1(\varOmega)$ onto $L^1(A)$ defined by
$P_A(f) = \chi_A f$.
\begin{definition}
    A subspace $X$ of $L^1(\varOmega)$ is said to be \emph{small}, if there is 
    no $A \in \varSigma$ of positive measure such that $P_A$ maps $X$ onto 
    $L^1(A)$.    
\end{definition}
If $X$ is a poor subspace of $L^1(\varOmega)$, then $X$ is small 
\cite{KadetsShepelskaWernerQuotientsDaugavetProperty}*{Corollary~6.7}.
The converse is valid too.
\begin{prop}
    \label{PropositionSmallPoor}
    If $X$ is a small subspace of $L^1(\varOmega)$, then $X$ is a poor
    subspace of $L^1(\varOmega)$.
    \begin{proof}
        Fix $A \in \varSigma$ with $\mu(A) > 0$ and $\varepsilon > 0$. By
        Proposition~\ref{PropositionPoorL1}, we have to find $f \in S_{L^\infty
        (\varOmega)}$ with $\supp f \subset A$ and $\lVert \restr{f}{X}\rVert
        \leq \varepsilon$.
        
        Since $X$ is small, the projection $P_A : L^1(\varOmega) \rightarrow
        L^1(A)$ does not map $X$ onto $L^1(A)$. By the (proof of the) open
        mapping theorem, $P_A[\frac 1 \varepsilon B_X]$ is nowhere dense in 
        $L^1(A)$. Pick $g \in B_{L^1(A)}$ with $g \notin 
        \overline{P_A[\frac 1 \varepsilon B_X]}$. The set 
        $\overline{P_A[\frac 1 \varepsilon B_X]}$ is absolutely convex and so 
        there exists by the Hahn-Banach theorem a function $f \in 
        S_{L^\infty(A)}$ with
        \begin{equation*}
            \sup \left\{ \left|\int_A fh\,d\mu\right| :
            h \in \frac 1 \varepsilon B_X \right\} \leq \Re \int_A fg
            \,d\mu.
        \end{equation*}
        Using this inequality, we get
        \begin{equation*}
            \norm{\restr{f}{X}} = \sup \left\{ \left| \int_A f h 
            \,d\mu \right| : h \in B_X \right\}
            \leq \varepsilon \Re \int_A fg \,d\mu \leq \varepsilon. \qedhere
        \end{equation*}
    \end{proof}
\end{prop}
An important tool in the study of small subspaces is the topology of convergence
in measure.

\mbox{}

\mbox{}
\begin{definition}
    \mbox{}
    \begin{enumerate}[\upshape (a)]
        \item A subspace $X$ of $L^1(\varOmega)$ is called \emph{nicely placed},
            if $B_X$ is closed with respect to convergence in measure.
        \item $\varLambda$ is said to be \emph{nicely 
            placed}, if $L^1_\varLambda(G)$ is a nicely placed subspace of
            $L^1(G)$.
        \item $\varLambda$ is said to be a \emph{Shapiro
            set}, if every subset of $\varLambda$ is nicely placed.
    \end{enumerate}
\end{definition}
G.~Godefroy coined these terms
\citelist{\cite{GodefroySousEspacesBienDisposes} \cite{GodefroyRieszSubsets}}
and showed that every Shapiro set is a Riesz set 
\cite{HarmandWernerMideale}*{Proposition~IV.4.5}. The natural numbers are a 
Shapiro set of $\mathbb Z$  \cite{HarmandWernerMideale}*{Example~IV.4.11} and
$\varLambda = \bigcup_{n = 0}^\infty \{k2^n: |k| \leq 2^n\}$ is a nicely placed
Riesz set which is not a Shapiro set 
\cite{HarmandWernerMideale}*{Example~IV.4.12}.
\begin{lem}
    \label{LemmaPoorInjection}
    Let $X$ be a nicely placed subspace of $L^1(G)$ and suppose that 
    there exists $A \in \mathcal B(G)$ with $m(A) > 0$ such that $P_A$ maps $X$
    onto $L^1(A)$, i.e., suppose that $X$ is not small. Then there exists a 
    continuous operator $T : L^1(A) \rightarrow X$ with 
    $j_A = P_AT$ where $j_A : L^1(A) \rightarrow L^1(G)$ is the natural
    injection.
    \begin{proof}
        This proof is a modification of a proof by G.~Godefroy, N.\,J.~Kalton, 
        and D.~Li \cite{GodefroyKaltonLiOperatorsBetweenSubspacesQuotients}*{
        Lemma~III.5}. We identify $X^{**}$ with $X^{\perp\perp} \subset
        L^1(G)^{**}$ and recall that A.\,V.~Buhvalov and 
        G.\,Ya.~Lozanovski{\u\i} showed that $P[B_{X^{\perp\perp}}] = B_X$, if 
        $X$ is nicely placed in $L^1(G)$
        \cite{HarmandWernerMideale}*{Theorem~IV.3.4}.
        
        Denote by $\mathcal N$ the directed set of open neighborhoods of $e_G$.
        (We turn $\mathcal N$ into a directed set by setting
        $V \leq W$ if and only if $V$ contains $W$.) Let $\mathcal U$ be a 
        ultrafilter on $\mathcal N$ which contains the filter base
        $\{ \{ W \in \mathcal N: V \leq W\} : V \in \mathcal N\}$.
        
        By the open mapping theorem, we can fix $M > 0$ with $B_{L^1(A)}
        \subset M P_A[B_X]$. For every $V \in \mathcal N$, use 
        Lemma~\ref{LemmaCovering} and choose
        disjoint Borel sets $B_{V,1},\dotsc,B_{V,N_V}$ with
        $A = \bigcup_{k = 1}^{N_V} B_{V,k}$ and $B_{V,k} - B_{V,k} \subset V$
        for $k = 1, \dotsc, N_V$. Picking $f_{V,k} \in MB_X$ with
        $P_A (f_{V,k}) = \frac{1}{m(B_{V,k})}\chi_{B_{V,k}}$ for 
        $k =1 , \dotsc, N_V$, we define $S_V : L^1(A) \rightarrow X$ by
        \begin{equation*}
            S_V (f) = \sum_{k = 1}^{N_V} \left(\int_{B_{V,k}} f \,dm \right) 
            f_{V,k} \quad (f \in L^1(A)).
        \end{equation*}
        Since the norm of every $S_V$ is bounded by $M$, we can define
        $S : L^1(A) \rightarrow X^{\perp\perp}$ by
        \begin{equation*}
            S (f) = \text{$w^*$-}\lim_{V,\mathcal U} S_V (f) 
            \quad (f \in L^1(A))
        \end{equation*}
        and set $T = PS$.
        
        Let us check that $j_A = P_AT$. Fix $f \in L^1(A)$. Since $C(G)$ is 
        dense in $L^1(G)$, we may assume that $f$ is the restriction to $A$
        of a continuous function. Let $(S_{\varphi(j)}(f))_{j \in J}$ be a 
        subnet of $(S_V (f))_{V \in \mathcal N}$ with
        $S(f) = \text{$w^*$-}\lim_j S_{\varphi(j)}(f)$. Since $f$ is
        uniformly continuous, it is easy to construct an
        increasing sequence $(j_n)_{n \in \mathbb N}$ in $J$ with
        \begin{equation}
            \label{LemmaPoorInjectionEqu1}
            \sup\left\{\norm{f - P_A S_{\varphi(j)}(f)}_\infty : j \geq j_n
            \right\} \longrightarrow 0.
        \end{equation}
        Furthermore, there exists a sequence $(g_n)_{n \in \mathbb N}$ in
        $L^1(G)$ that converges $m$-almost everywhere to $PS(f)$ with 
        $g_n \in \co\{ S_{\varphi(j)}(f):j \geq j_n \}$ for all 
        $n \in \mathbb N$ \cite{HarmandWernerMideale}*{Lemma~IV.3.1}. Hence we 
        have by (\ref{LemmaPoorInjectionEqu1}) that for $m$-almost all $x \in A$ 
        \begin{equation*}
            T(f)(x) = PS(f)(x) = \textstyle \lim_n g_n(x) = f(x) 
        \end{equation*}
        and therefore $j_A = P_AT$.
    \end{proof}
\end{lem}
\begin{theo}
    \label{TheoremNicelyPlacedRieszPoor}
    If $\varLambda$ is a nicely placed Riesz set, then $L^1_\varLambda(G)$ is a
    small subspace of $L^1(G)$.
    \begin{proof}
        Assume that $\varLambda$ is a nicely placed Riesz set such that 
        $L^1_\varLambda(G)$ is not a small subspace of $L^1(G)$. 
        
        Since $L^1_\varLambda(G)$ is not small, there exists a Borel set $A$ of 
        positive measure such that $P_A$ maps $L^1_\varLambda(G)$ onto $L^1(A)$. 
        Using Lemma~\ref{LemmaPoorInjection}, we find $T : L^1(A) \rightarrow
        L^1_\varLambda(G)$ with $j_A = P_AT$. This operator is an isomorphism
        onto its image and $L^1_\varLambda(G)$ contains a copy of $L^1(A)$. So
        $L^1_\varLambda(G)$ fails the Radon-Nikod\'ym property. But this
        contradicts our assumption because $\varLambda$ is a Riesz set if and 
        only if $L^1_\varLambda(G)$ has the Radon-Nikod\'ym property
        \cite{LustEnsemblesRosenthalRiesz}*{Th\'eor\`eme~2}.
    \end{proof}
\end{theo}
\begin{cor}
    If $\varLambda$ is a Shapiro set, then $L^1_\varLambda(G)$ is a poor 
    subspace of $L^1(G)$.
\end{cor}
Theorem~\ref{TheoremNicelyPlacedRieszPoor} can be strengthened, if 
$G$ is metrizable. Let $\varLambda$ be nicely placed. Then $L^1_\varLambda(G)$ 
s a poor subspace of $L^1(G)$ if and only if $\varLambda$ is a semi-Riesz set
\cite{GodefroyKaltonLiOperatorsBetweenSubspacesQuotients}*{Proposition~III.10}

\section*{Acknowledgements}

This is part of the author's Ph.D. thesis, written under the supervision of
D.~Werner at the Freie~Universit\"at~Berlin.

\end{document}